\documentclass{article}
\usepackage{hyperref}
\usepackage{amsmath,amsthm,amssymb,enumerate,latexsym,verbatim,amscd}
\usepackage{pstricks,pst-node,pst-tree}

\newenvironment{enumeratei}
{\begin{enumerate}[\upshape (i)]}
{\end{enumerate}}

\newtheorem{thm}{Theorem}[section]
\newtheorem{prop}[thm]{Proposition}
\newtheorem{lem}[thm]{Lemma}
\newtheorem{cor}[thm]{Corollary}

\theoremstyle{remark}
\newtheorem{rmk}[thm]{Remark}

\newcommand{\nc}{\newcommand}
\nc{\hb}{\mathbb}
\nc{\mbf}{\mathbf}
\nc{\DMO}{\DeclareMathOperator}
\nc{\MyNode}[1]{\TR{\makebox{(#1)}}}

\nc{\cA}{{\mathcal A}} \nc{\cB}{{\mathcal B}} \nc{\cC}{{\mathcal C}}
\nc{\cD}{{\mathcal D}} \nc{\cE}{{\mathcal E}} \nc{\cF}{{\mathcal F}}
\nc{\cG}{{\mathcal G}} \nc{\cH}{{\mathcal H}} \nc{\cI}{{\mathcal I}}
\nc{\cJ}{{\mathcal J}} \nc{\cK}{{\mathcal K}} \nc{\cL}{{\mathcal L}}
\nc{\cM}{{\mathcal M}} \nc{\cN}{{\mathcal N}} \nc{\cO}{{\mathcal O}}
\nc{\cP}{{\mathcal P}} \nc{\cQ}{{\mathcal Q}} \nc{\cR}{{\mathcal R}}
\nc{\cS}{{\mathcal S}} \nc{\cT}{{\mathcal T}} \nc{\cU}{{\mathcal U}}
\nc{\cV}{{\mathcal V}} \nc{\cW}{{\mathcal W}} \nc{\cX}{{\mathcal X}}
\nc{\cY}{{\mathcal Y}} \nc{\cZ}{{\mathcal Z}}

\nc{\Aa}{{\hb A}} \nc{\Cc}{{\hb C}} \nc{\Gg}{{\hb G}}
\nc{\Nn}{{\hb N}} \nc{\Pp}{{\hb P}} \nc{\Qq}{{\hb Q}}
\nc{\Rr}{{\hb R}} \nc{\Zz}{{\hb Z}}

\nc{\sB}{{\mathsf B}} \nc{\sE}{{\mathsf E}}
\nc{\sW}{{\mathsf W}} \nc{\sX}{{\mathsf X}}
\nc{\sY}{{\mathsf Y}}

\nc{\tuple}[2]{{#1},\ldots,{#2}} \nc{\ptu}[2]{{#1}:\ldots:{#2}}
\nc{\pbrac}[1]{\langle {#1} \tongle}

\nc{\maps}[3]{{#1}\colon {#2} \to {#3}}
\nc{\map}[2]{{#1}\to {#2}}
\nc{\res}[2]{{#1} |_{#2}}
\nc{\imbed}{\hookrightarrow}

\nc{\set}[1]{\{#1\}}

\nc{\bn}{\binom}

\nc{\bds}{\boldsymbol}

\DMO{\h}{ht}
\DMO*{\pgl}{{\sf PGL}}
\DMO*{\psl}{{\sf PSL}}
\DMO*{\gl}{{\sf GL}}
\nc{\norm}[1]{\|#1\|}

\title{Diophantine Equations of Matching Games I} 

\author{Chun Yin Hui\\
  Department of Mathematics,\\
  Indiana University Bloomington\\
  and \\
  Wai Yan Pong\\
  Department of Mathematics,\\ California University Dominguez Hills}

\begin{document}
\maketitle
\begin{abstract}
  We solve a family of quadratic Diophantine equations associated to a
  simple kind of games. We show that the ternary case, in many ways,
  is the most interesting and the least arbitrary member of the
  family.
\end{abstract}

\section{The Matching Games} 
\label{s:game} 

An {\em $(n,d)$-matching game} ($n,d \ge 2$) is a game in which the
player draws $d$ balls from a bag of balls of $n$ different
colors. The player wins if and only if the balls drawn are all of the
same color. A game is {\em non-trivial} if there are at least $d$
balls in the bag. It is {\em faithful} if there are balls in each of
the $n$ colors. A game is {\em fair} if the player has an equal chance
of winning or losing the game. In this article, we only study the
$(n,2)$-matching games or simply the {\em $n$-color games}, leaving
the study of the higher $d$ case to~\cite{dioII}.

An $n$-tuple $(\tuple{a_1}{a_n})$ where $a_i$ is the number of the
$i$-th color balls in the bag represents an $n$-color game. For $m \le
n$, an $m$-color game $(\tuple{a_1}{a_m})$ can be regarded as the
$n$-color game $(\tuple{a_1}{a_m},\tuple{0}{0})$. The only trivial
$n$-color fair games, are the {\em zero game} $(0, \ldots, 0)$ and, up
to permutation, the game $(0,0, \ldots, 1)$.

By considering the number of ways for the player to win the game, one
sees that the $n$-color fair games are exactly the non-negative
integral solutions of
\begin{equation*}
  \binom{\sum_{i=1}^n x_i}{2} = 2\left(\sum_{i=1}^n
    \binom{x_i}{2}\right)
\end{equation*}
or equivalently,
\begin{equation}
  \label{eq:fair}
  F_n(x_1,\ldots, x_n) := \left(\sum\nolimits_{i =1}^n x_i\right)^2 -
  \sum\nolimits_{i=1}^n x_i - 4\sum\nolimits_{i \neq j} x_ix_j =0
\end{equation}

The paper will be organized in the following way: We give a brief
treatment of the 2-color games in Section~\ref{s:n=2}. The results
there illustrate the kind of questions that we try to answer in the
general case. In Section~\ref{s:n>=3}, we give a ``parametric''
solution to Equation~\eqref{eq:fair}. It is unclear, however, from
this method which choice of the parameters will yield the fair
games. We tackle this problem in Section~\ref{s:graph} by giving a
graph structure to the solutions.  We show that the components of this
graph are trees and give an algorithm for finding their roots. This
yields all solutions recursively. Furthermore, we characterize the
components containing the fair games. For $n=3$, we show that the
graph consists of two trees with the nontrivial 3-color fair games
forming a full binary tree. We then study what are the possible
coordinates of fair games in Section~\ref{s:C_n}. In
Section~\ref{s:asy}, we establish some partial results concerning the
asymptotic behavior of 3-color fair games. We conclude the article
with some odds and ends of our study in Section~\ref{s:misc}.

The following conventions will be used throughout this article:
\begin{itemize}
\item All variables and unknowns range over the integers unless
  otherwise stated.

\item The cardinality of a set $A$ is denoted by $|A|$. 

\item For $\bds{a} \in \Zz^m$, the (Euclidean) norm of $\bds{a}$ is
  denoted by $\norm{\bds{a}}$. For $A \subseteq \Zz^m$ and $k \ge 0$,
  $A(k)$ denotes the set of elements of $A$ with norm at most $k$. We
  define the {\em height} of $\bds{a}$ so be
  $\left|1+\sum{a_i}\right|$.
  
\item For any integer $d$ and $\bds{a}, \bds{b} \in \Zz^m$, we say
  that $\bds{a}$ and $\bds{b}$ are congruent modulo $d$, written as
  $\bds{a} \equiv \bds{b} \mod d$, if $a_i \equiv b_i \mod d$ for all
  $1 \le i \le m$.
  
\item Denote by $\cS_n$ and $\cF_n$ the set of integral and
  non-negative integral solutions, identified up to permutations, of
  Equation~\eqref{eq:fair} respectively. Elements of $\cF_n$ are the
  $n$-color fair games. We often use an increasing (or decreasing)
  tuple to represent an element of $\cS_n$. Denote by $\cC_n$ the set
  of coordinates of $\cF_n$.

\item For any $n$-tuple $\bds{x} = (x_1,\ldots, x_n)$ and $I$ a subset
  of the indices, we write $\bds{x}_I$ for the tuple obtained from
  $\bds{x}$ by omitting the variables indexed by the elements of
  $I$. We write $\bds{x}_i$ for $\bds{x}_{\{i\}}$ and $\bds{x}_{ij}$
  for $\bds{x}_{\{i,j\}}$, etc.

\item Let $s(\bds{x})$ and $p(\bds{x})$ be the symmetric polynomials
  of degree 1 and 2, respectively, i.e.
  \[
  s(\bds{x}) = \sum\nolimits_{i =1}^n x_i, \qquad p(\bds{x}) =
  \sum\nolimits_{1 \le i < j \le n} x_ix_j
  \]
  We often omit writing out the variables explicitly, so we write
  $s_i$ for $s(\bds{x}_i)$, $s_{ij}$ for $s(\bds{x}_{ij})$, etc. We
  understood $s \equiv 0$ on zero variables and $p \equiv 0$ on either
  0 or 1 variable.
\end{itemize}

Thanks go to Jackie Barab from whom the second author first learned
about the 2-color games~\footnote{They are used in teaching 3rd and
  4th graders in California about probability}. We thanks Bjorn Poonen
for referring~\cite{gs} to us. We would also like to thank Thomas
Rohwer and Zeev Rudnick for bringing~\cite{ntd} and~\cite{lorentzian},
respectively, to our attention. Finally, we thank Michael Larsen for
reading a draft of this article and giving us several valuable
comments.

\section{The 2-color games}
\label{s:n=2}
As a warm-up, we analyze the 2-color games first. In this case,
Equation~\eqref{eq:fair} becomes $(x_1 - x_2)^2-(x_1 + x_2)=0$ and is
easy to solve: let $m=x_2-x_1 \ge 0$, and then $x_1 = m(m-1)/2$ and so
$x_2 = (m+1)m/2$. This shows that
\begin{thm}
  \label{t:n=2}
  The $2$-color fair games are pairs of consecutive triangular
  numbers. In particular, $\cC_2$ is the set of triangular numbers.
\end{thm}
Using Theorem~\ref{t:n=2}, a few simple computations tell us the number
of fair 2-color games of norm bounded by a given number.
\begin{cor}
  \label{c:asymn=2}
  For $k \ge 0$,
  \begin{enumerate}
  \item \label{i:norm-n=2} $|\cF_2(k)|=[\sqrt{r(k)}]+1$. Hence
    $|\cF_2(k)|$ is asymptotic to $2^{1/4}\sqrt{k}$.
    
  \item \label{i:C2} $|\cC_2(k)|=[r(\sqrt{k})]$. Hence $|\cC_2(k)|$ is
    asymptotic to $\sqrt{2k}$.
  \end{enumerate}
\end{cor}
Here $r(k) = (-1 + \sqrt{1+8k^2})/2$ and $[x]$ is the largest integer
$\le x$.

\section{Solving Equation~\eqref{eq:fair}} 
\label{s:n>=3} 
It would be nice to know in advance that Equation~\eqref{eq:fair} is
solvable. The following simple observation tells us just that.
\begin{thm}
  \label{t:solvable}
  There are infinitely many faithful $n$-color fair games.
\end{thm}
\begin{proof}
  Regarding the polynomial $F_n$ in~\eqref{eq:fair} as a quadratic in
  $x_k$, we have
  \begin{equation}
    \label{eq:x_k}
    F_n(\bds{x}) = x_k^2 - (2s_k +1)x_k + F_{n-1}(\bds{x}_k).
  \end{equation}
  Thus if $(a_1, \ldots, a_k, \ldots, a_n)$ is a solution then so is
  $(a_1, \ldots, b_k, \ldots, a_n)$ where $b_k = 2s_k(\bds{a}) +1
  -a_k$. In particular, if $(a_1, \ldots, a_{n-1})$ is a fair game,
  then so are $(a_1, \ldots, a_{n-1}, 0)$ and $(a_1, \ldots, a_{n-1},
  1+2\sum_{i < n} a_i)$. The latter game is faithful if $(a_1, \ldots,
  a_{n-1})$ is. Since there are infinitely many faithful 2-color fair
  games (Theorem~\ref{t:n=2}), the theorem follows by induction on
  $n$.
\end{proof}
Since the two roots of~\eqref{eq:x_k} sum to $2s_k+1$, they can be
expressed as $s_k + m +1$ and $s_k-m$ for some $m \ge 0$. Thus solving
\begin{equation}
  \label{eq:disc}
  \begin{split}
    (s_k +m +1)(s_k -m) &= F_{n-1}(\bds{x}_k) =s_k^2 - s_k -4p_k \\
    2s_k + 4p_k &= m^2 + m
  \end{split}
\end{equation}
will solve~\eqref{eq:fair} and vice versa. After adding $1+4s_{ijk}^2 +
2s_{ijk} - 4p_{ijk}$ ($i,j,k$ pairwise distinct) to both sides
of~\eqref{eq:disc}, the left-hand side factorizes:
\begin{equation}
  \label{eq:ff}
  \begin{split}
    (2x_i + 2s_{ijk} +1)&(2x_j +2s_{ijk} +1) = m^2 + m +1 + 4s_{ijk}^2
    + 2s_{ijk} - 4p_{ijk} \\
    & = m^2 + m + 1 + 2(s_{ijk}^2 + s_{ijk} + \norm{\bds{x}_{ijk}}^2)
  \end{split}
\end{equation}
Equation~\eqref{eq:ff} gives us a way to solve
Equation~\eqref{eq:fair}: Denote by $J(\bds{x}_{ijk},m)$ the
right-hand side of~\eqref{eq:ff}. Choose $\bds{a} \in \Zz^{n-3}$ and
$m \ge 0$ arbitrarily\footnote{When $n=3$, we only need to choose
  $m$.}. According to~\eqref{eq:ff}, we can solve for $x_i$ and $x_j$
by factorizing the odd number $J(\bds{a},m)$ into a product two odd
numbers. By~\eqref{eq:x_k}, we can then solve for $x_k$ and hence a
solution of~\eqref{eq:fair}. Moreover, it is clear from the discussion
above that any solution of~\eqref{eq:fair} arises from such a
factorization. For the record, we have
\begin{thm}
  \label{t:S_n}
  Fix $n \ge 3$. For any $\bds{a} \in \Zz^{n-3}$, $m \ge 0$ and $0 \le
  b \le c$ such that $J(\bds{a},m)=(2b+1)(2c+1)$, the following are
  solutions to Equation~\eqref{eq:fair}:
  \begin{align*}
    & (b-s(\bds{a}),\ c-s(\bds{a}),\ b+c+1-s(\bds{a}) +m,\ \bds{a}) \\
    & (b-s(\bds{a}),\ c-s(\bds{a}),\ b+c-s(\bds{a})-m,\ \bds{a})\\
    & (-(c+1)-s(\bds{a}),\ -(b+1)-s(\bds{a}),\ -(b+c+
    s(\bds{a})+1) +m,\ \bds{a}) \\
    & (-(c+1)-s(\bds{a}),\ -(b+1)-s(\bds{a}),\ -(b +c +
    s(\bds{a})+2) -m,\ \bds{a})
  \end{align*}
  Moreover, up to a permutation every solution of
  Equation~\eqref{eq:fair} is in one of these forms.
\end{thm}
There is a less tricky way to derive Equation~\eqref{eq:ff}. We give
the idea here but leave the details to the
reader. Equation~\eqref{eq:disc} can be viewed as a curve on the
$x_ix_j$-plane ($i,j,k$ pairwise distinct). One can express the curve
having an integral point by expressing that the corresponding
quadratic in $x_i$ is solvable in terms of $\bds{x}_{ijk}$ and $d :=
x_j -x_i$. The expression $1+4s_{ijk}^2 + 2s_{ijk} - 4p_{ijk}$ then
flows out naturally.

Even though the method above solves Equation~\eqref{eq:fair}, it is
unclear which choice of the parameters will produce fair games. For
example, $J(2,3) =33$ does not produce any 4-color fair game. We will
take up this issue in the next section.

\section{Solutions as a graph}
\label{s:graph}
Starting from a solution $(\tuple{a_1}{a_n})$ of
Equation~\eqref{eq:fair}, we obtain another one by replacing $a_k$
with $b_k := 2\sum_{i \neq k} a_i + 1 - a_k$ (see
Theorem~\ref{t:solvable}). This suggests that we can view $\cS_n$ as a
graph by putting an edge between two elements of $\cS_n$ if they
differ at only one coordinate\footnote{Incidentally, this is the same
  graph structure that was put on the solutions on the Markoff
  Equation~\cite{cus,mar}.}. An immediate question would be: can we
generate every fair game from some fixed game, say the zero game?  In
other words, is $\cF_n$ connected as a graph?  We have seen that
$\cF_2$ is connected (Theorem~\ref{t:n=2}) and we will show that the
same is true for $\cF_3$. However, $\cF_n$ fails to be connected for
$n \ge 4$.

Let us begin with a crucial observation. For any $\bds{a} \in \cS_n$
and any three pairwise distinct indices $i,j,k$, according
to~\eqref{eq:ff}, for some $m \ge 0$,
\[
(2s_{jk}(\bds{a}) +1)(2s_{ik}(\bds{a})+1) = 2(s_{ijk}(\bds{a})^2 +
s_{ijk}(\bds{a}) + \norm{\bds{a}_{ijk}}^2)+m^2+m+1.
\]
Since the right-hand-side is always positive, we conclude that
\begin{prop}
  \label{p:sum n-2}
  For any $\bds{a} \in \cS_n$, the numbers $2s_{ij}(\bds{a})+1$
  \upshape{(}$1 \le i < j \le n$\upshape{)} all have the same sign. In
  particular, the coordinates of an element of $\cS_3$ are either all
  non-negative or all negative.
\end{prop}
We define the {\em sign} of $\bds{a}$ as the common sign of the
$2s_{ij}(\bds{a}) +1$ ($1 \le i < j \le n$). Note that it is the same of
the sign of $s(\bds{a})+1$ since $\sum_{i <j} s_{ij} =
\binom{n-1}{2}s$. Let $\cS_n^+$ and $\cS_n^-$ be the sets of positive
and negative elements of $\cS_n$, respectively. Since any two neighbors
in $\cS_n$ share $n-1$ coordinates, they must have the same sign,
therefore
\begin{prop}
  \label{p:nbd-sign}
  $\cS_n^+$ and $\cS_n^-$ are disjoint union of components of $\cS_n$.
\end{prop}
Our next result shows how height varies among neighbors.
\begin{prop}
  \label{p:height}
  At most one neighbor of any vertex of $\cS_n$ can have a smaller
  height. Moreover, any two neighbors must have different height.
\end{prop}
\begin{proof}
  Fix any $\bds{a} \in \cS_n$. Let $\bds{b}_k = (a_1, \ldots, b_k,
  \ldots, a_n)$ where $b_k = 2s_k(\bds{a}) +1 - a_k$ $(1 \le k \le n)$
  be its neighbors. Rearranging the coordinates if necessary, we
  assume $a_1 \le a_2 \le \cdots \le a_n$.
  \begin{par}
    \noindent {Case 1: $\bds{a} \in \cS_n^+$.} Then for $k \neq n$,
    \[
    b_k=2s_{kn}(\bds{a}) + (a_n - a_k) +a_n +1 > a_n \ge a_k
    \]
    so $s(\bds{b}_k) > s(\bds{a}) \ge 0$. \medskip
    
    \noindent {Case 2: $\bds{a} \in \cS_n^-$.} Then for $k \neq 1$,
    \[
    b_k = 2s_{1k}(\bds{a}) + (a_1 -a_k) + a_1 + 1 < a_1 \le a_k
    \]
    so $s(\bds{b}_k) < s(\bds{a}) < 0$.
  \end{par}

  \noindent This completes the proof of the first statement since in
  both cases we have $\h(\bds{b}_k) > \h(\bds{a})$ for all but perhaps
  one $k$. The second statement follows readily from the fact that each
  $a_k+b_k$ is an odd number.
\end{proof}
We say that a vertex of $\cS_n$ is a {\em root} if all its neighbors
have a greater height. We would like to point out that replacing
height by norm in the definition of root will yield the same concept
since Equation~\eqref{eq:fair} can be rewritten as
  \begin{equation}
    \label{eq:height-norm}
    \binom{s+1}{2} = \norm{\bds{x}}^2.
  \end{equation}
\begin{thm}
  \label{t:tree}
  The components of $\cS_n$ are rooted trees.
\end{thm}
\begin{proof}
  Proposition~\ref{p:height} implies that for any subgraph $H$ of
  $\cS_n$, a vertex of maximal height in $H$ cannot have two neighbors
  in $H$. This shows that $\cS_n$ must be acyclic. Moreover, every
  component of $\cS_n$ has a unique vertex of minimal height. If not,
  take a path with two vertices of minimal height as endpoints. Since
  neighbors in $\cS_n$ have different heights, the path has length at
  least 2 but then a vertex of maximal height in the path will have
  two neighbors, a contradiction.
\end{proof}
\begin{figure}
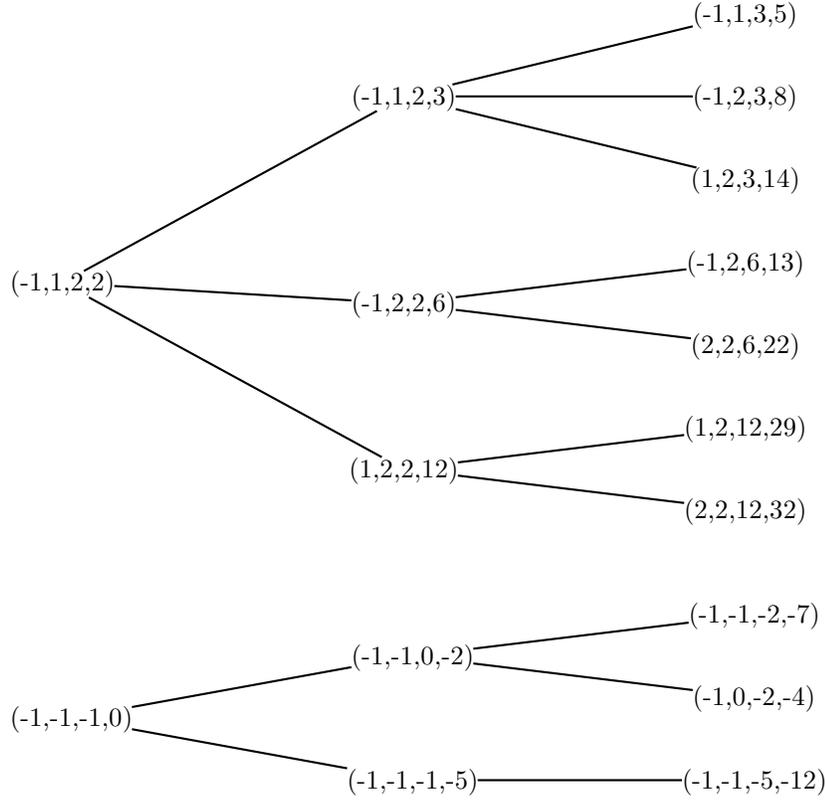

  \pstree[treemode=R, levelsep=30ex]{\MyNode{-1,1,2,2}}
  {\pstree{\MyNode{-1,1,2,3}} 
          {\MyNode{-1,1,3,5} \MyNode{-1,2,3,8} \MyNode{1,2,3,14}}
   \pstree{\MyNode{-1,2,2,6}} {\MyNode{-1,2,6,13} \MyNode{2,2,6,22}}
   \pstree{\MyNode{1,2,2,12}} {\MyNode{1,2,12,29} \MyNode{2,2,12,32}}}
  
  \vspace{1cm}

  \pstree[treemode=R, levelsep=30ex]{\MyNode{-1,-1,-1,0}}
  {\pstree{\MyNode{-1,-1,0,-2}} {\MyNode{-1,-1,-2,-7}
      \MyNode{-1,0,-2,-4}} 
   \pstree{\MyNode{-1,-1,-1,-5}}
    {\MyNode{-1,-1,-5,-12}}}
  \caption{A positive tree and a negative tree in $\cS_4$}
\end{figure}
Fair games are positive solutions of~\eqref{eq:fair} and yet a
positive solution, for example (-1,1,2,2), need not even represents a
game. However, for any $\bds{a} \in \cS_n^{+}$, a neighbor of
$\bds{a}$ greater in height will have the different coordinate
non-negative (see the proof of Proposition~\ref{p:height}). Thus by
going up in height along any branch, we see that
\begin{prop}
  \label{p:positive}
  Every component of $\cS_n^{+}$ contains fair games.
\end{prop}
\noindent Similarly, by going down in height, we see that every component of
$\cS_n^{-}$ contains solutions with all negative coordinates.

Our next goal is to locate the roots of $\cS_n$. Once this is achieved,
we will have an effective way of generating all fair games since each of
them is connected to some positive root.
\begin{prop}
  \label{p:root_bound}
  Suppose $\bds{r} \in \cS_n^{+}$ {\upshape (}$\cS_n^{-}$,
  resp.{\upshape )} is a root and $i,j,k$ pairwise distinct where $k$
  is the index of a maximal {\upshape (}minimal, resp.{\upshape )}
  coordinate of $\bds{r}$.  Then $\bds{r}$ is obtained from a
  factorization of $J(\bds{r}_{ijk},m)$ for some $0 \le m \le
  B(\bds{r}_{ijk})$ where $B(\bds{r}_{ijk})$ is an explicit bound give
  in terms of $\bds{r}_{ijk}$.
\end{prop}
\begin{proof}
  We argue for $S_n^+$ only. The proof for $S_n^-$ is similar. According
  to~\eqref{eq:ff}, for some $m \ge 0$,
  \[
  (2s_{ik}+1)(2s_{jk}+1) = 2(s_{ijk}^2 + s_{ijk}
  +\norm{\bds{r}_{ijk}}^2) + m^2 +m +1.
  \]
  Since $\bds{r}$ is a root, $r_k$ is the smaller root of
  Equation~\eqref{eq:x_k}, i.e. $r_k = s_k-m$. And since $r_k \ge
  r_{\ell}$ ($\ell \neq k$), so $s_{\ell k} \ge s_{\ell k} +r_{\ell}
  -r_k = m$. Thus
  \begin{align*}
    (2m+1)^2 &\le 2(s_{ijk}^2 + s_{ijk} + \norm{\bds{r}_{ijk}}^2) +
    m^2
    + m +1 \\
    3m^2 +3m &\le 2(s_{ijk}^2 + s_{ijk} + \norm{\bds{r}_{ijk}}^2).
  \end{align*}
  That means $0 \le m \le B(\bds{r}_{ijk})$ where $B(\bds{r}_{ijk})$
  expresses the larger root of the quadratic $3x^2 + 3x - 2(s_{ijk}^2
  + s_{ijk} + \norm{\bds{r}_{ijk}}^2)$.
\end{proof}
Let us summarize how to find the roots of $\cS_n$: for each $\bds{a} \in
\Zz^{n-3}$, we compute the finite set consisting of those solutions
given by the factorizations of $J(\bds{a},m)$ where $0 \le m \le
B(\bds{a})$. We then check which element in this finite set is a
root. While Proposition~\ref{p:root_bound} guarantees that every root of
$\cS_n$ can be found this way, our next result shows that we do have to
check for every $\bds{a} \in \Zz^3$.
\begin{prop}
  \label{p:roots}
  Every $(n-3)$-tuple of integers can be extended to a root in
  $\cS_n$. More precisely, for any $\bds{a} \in \Zz^{n-3}$, $n$-tuples
  \begin{align*}
    & \bds{r}_{+}:=(s(\bds{a})^2 + \norm{\bds{a}}^2,\ s(\bds{a})^2 +
    \norm{\bds{a}}^2,\
    -s(\bds{a}),\ \bds{a}) \quad \text{and}\\
    & \bds{r}_{-}:=(-(s(\bds{a})+1)^2 - \norm{\bds{a}}^2,\
    -(s(\bds{a})+1)^2 - \norm{\bds{a}}^2,\ -(s(\bds{a}) +1),\ \bds{a})
  \end{align*}
  are a positive and a negative root of $\cS_n$, respectively.
\end{prop}
\begin{proof}
  First by Theorem~\ref{t:S_n}, they are solutions corresponding to
  the trivial factorization of $J(\bds{a},0)$ (in the notation there,
  $b = s(\bds{a})^2 + s(\bds{a}) + \norm{\bds{a}}^2$ and $c
  =0$.). Clearly, $\bds{r}_{+}$ is positive while $\bds{r}_{-}$ is
  negative. The neighbor of $\bds{r}_{+}$ obtained by varying its
  largest coordinate $s(\bds{a})^2 + \norm{\bds{a}}^2$ has an even
  larger coordinate, namely $s(\bds{a})^2 + \norm{\bds{a}}^2 +
  1$. Thus $\bds{r}_{+}$ is indeed a root (see the proof of
  Theorem~\ref{p:height}). A similar argument show that $\bds{r}_{-}$
  is a root as well.
\end{proof}
Since each component of $\cS_n$ has exactly one root, an immediate
consequence of Proposition~\ref{p:roots} is that
\begin{thm}
  \label{t:inf-many-comps} For $n \ge 4$, $\cS_n^{+}, \cS_n^{-}$ each
  has infinitely many components.
\end{thm}
On the contrary, by examining the proof of
Proposition~\ref{p:root_bound}, one readily checks that $(0,0,0)$ is the
only root of $\cS_3^{+}$. By Theorem~\ref{t:S_n}, the map $\bds{a}
\mapsto -(\bds{a} + \bds{1})$ where $\bds{1} = (1,1,1)$ is a graph
isomorphism between $\cS_3^{+}$ and $\cS_3^{-}$. Moreover, every element
of $\cS_3^{+}$ is actually a fair game according to
Proposition~\ref{p:sum n-2}. Thus,
\begin{thm}
  \label{t:S_3-comp}
  $\cS_{3}^{+}$ and $\cS_3^{-}$ are the two components of
  $\cS_3$. Moreover, $\cS_3^{+} = \cF_3$.
\end{thm}
A straight-forward computation shows that every vertex of $S_3$ with
distinct coordinates has two distinct children (i.e. neighbors with a
bigger norm). Moreover, each of its children also has distinct
coordinates. Hence,
\begin{thm}
  \label{t:tree S_3}
  The non-trivial 3-color fair games form an infinite full binary tree
  with $(0,1,3)$ as root.
\end{thm}
\begin{figure}
  \centering \pstree[treemode=R, levelsep=12ex] {\MyNode{0,0,0}}
  {\pstree{\MyNode{0,0,1}} {\pstree{\MyNode{0,1,3}}
      {\pstree{\MyNode{0,3,6}} {\pstree{\MyNode{0,6,10}}
          {\MyNode{0,10,15} \MyNode{6,10,33}} \pstree{\MyNode{3,6,19}}
          {\MyNode{3,19,39} \MyNode{6,19,48}} }
        
        \pstree{\MyNode{1,3,9}} {\pstree{\MyNode{1,9,18}}
          {\MyNode{1,18,30} \MyNode{9,18,54}} \pstree{\MyNode{3,9,24}}
          {\MyNode{3,24,46} \MyNode{9,24,64}} } } } }
\caption{Part of $\cF_3$}
\end{figure}

\section{The set $\cC_n$} 
\label{s:C_n} 
\begin{prop}
  \label{p:C_n}
  For $n \ge 4$, $\cC_n$ is the set of non-negative integers.
\end{prop}
\begin{proof}
  For any $a \ge 0$, let $\bds{a}$ be the $(n-3)$-tuple with all
  coordinates equal $a$. Then the child $(4(s(\bds{a})^2 +\norm{a}^2)
  + 3s(\bds{a}) +1, s(\bds{a})^2 +\norm{\bds{a}}^2, s(\bds{a})^2
  +\norm{\bds{a}}^2, \bds{a})$ of the positive root $\bds{r}_{+}$ in
  Proposition~\ref{p:roots} is a fair game with $a$ as a
  coordinate. Incidentally, this also shows that for $n \ge 4$, every
  natural number is a coordinate of some faithful $n$-color fair game.
\end{proof}
This leaves us only $\cC_3$ to study. It turns out that our analysis
of $\cC_3$ will yield another way of finding the 3-color fair games
(Theorems~\ref{t:3game} and~\ref{t:partition} ). First, note that
$\cC_3$ is the set of $c \ge 0$ such that the curve defined by
\begin{equation}
  \label{eq:xy}
  (x_1-x_2)^2 -(2c+1)(x_1+x_2) +c(c-1) = 0
\end{equation}
has a non-negative point. Arguing mod $2$, one sees that any integral
point on the parabola
\begin{equation}
  \label{eq:uv}
  u^2 -(2c+1)v +c(c-1) =0
\end{equation}
must have coordinates with the same parity. Thus, the transformation
$u=x_1-x_2,\ v=x_1+x_2$ is a 1-to-1 correspondence between the
integral points of these curves. Moreover, those $(x_1,x_2)$'s with
$x_1,x_2 \ge 0$ correspond to the $(u,v)$'s with $u \le v$. However,
the inequality is automatic:
\begin{prop}
  \label{p:uvsol}
  Solutions of Equation~\eqref{eq:uv} are of the form
  \begin{equation*}
    (u,v) = \left(u,\quad \frac{u^2+c(c-1)}{2c+1}\right)
  \end{equation*}
  where $u^2 \equiv -c(c-1) \mod (2c+1)$. In particular,~\eqref{eq:uv}
  is solvable if and only if $-c(c-1)$ is a quadratic residue mod
  $(2c+1)$. Moreover, $|u| \le v$ for every integral solution $(u,v)$.
\end{prop}
\begin{proof}
  The first statement is clear by considering Equation~\eqref{eq:uv}
  mod $(2c+1)$.  Note that for any $u$,
  \begin{equation*}
    -\frac{8c+1}{4} \le u^2 \pm (2c+1)u +c(c-1).
  \end{equation*}
  Since $c \ge 0$,
  \begin{equation*}
    \pm u -\frac{8c+1}{8c+4} \le \frac{u^2 +c(c-1)}{2c+1}.
  \end{equation*}
  So in particular
  \begin{equation*}
    |u| \le \frac{u^2+c(c-1)}{2c+1}
  \end{equation*}
  if both sides are integers. Therefore, $|u| \le v$ for any integral
  solution of~\eqref{eq:uv}.
\end{proof}
Solving for $x_1,x_2$ in terms of $u,v$ yields a parametrization of
$\cF_3$.
\begin{thm}
  \label{t:3game}
  The 3-color fair games are of the form
  \begin{equation}
    \label{eq:3para}
    \left(\frac{u^2+(2c+1)u +c(c-1)}{2(2c+1)},\quad \frac{u^2
    -(2c+1)u+c(c-1)}{2(2c+1)},\quad  c \right)
  \end{equation}
  where $c \in \cC_3$ and $u^2 \equiv -c(c-1) \mod (2c+1)$.
\end{thm}

For $i=-1,0,1$, let $P_{i}$ be the set of primes that are congruent to
$i$ mod $3$. Let $P_{i}^{\ge 0}$ be the set of natural numbers whose
prime factors are all in $P_{i}$. With this notation, we have
\begin{thm}
  \label{t:partition}
  $\cC_3 = \set{c \colon c \ge 0, \ 2c+1 \in P_{1}^{\ge 0} \cup
    3P_{1}^{\ge 0}}$.
\end{thm}
\begin{proof}
  By Proposition~\ref{p:uvsol} and the discussion preceding it,
  $\cC_3$ is the set of all $c$ such that $-c(c-1)$ is a quadratic
  residue mod $2c +1$. Since $2c+1$ is odd, $-c(c-1)$ and $-4c(c-1)$
  are either both squares or both non-squares mod $2c+1$.  Note that
  $-4c(c-1)$ is congruent to $-3$ mod $2c+1$.  Therefore $c \in \cC_3$
  if and only if $-3$ is a square mod $2c+1$.  This condition is
  equivalent to: for every prime $p$, $-3$ is a square mod $p^{v_p}$
  where $v_p$ is the exponent of $p$ in $2c+1$. That means $v_3=0$ or
  $1$ and $-3$ is a quadratic residue mod $p$ for $p \neq 3$. The law
  of quadratic reciprocity then yields
  \begin{equation*}
  \left(\frac{-3}{p}\right)=\left(\frac{-1}{p}\right)\left(\frac{3}{p}\right)
  = (-1)^{\frac{1}{2}(p-1)}\left(\frac{3}{p}\right)=\left(\frac{p}{3}\right).
  \end{equation*}
  So the last condition is equivalent to $p \equiv 1 \mod 3$. This
  completes the proof.
\end{proof}

With the binary tree structure of $\cF_3$ in mind, it is tempting to
establish some relation between the equation for fair 3-color games
with the Markoff equation (see~\cite{cus} and~\cite{mar} for more
information on the Markoff equation). Unfortunately, we could not find
any, however the naive analog of the Markoff conjecture (that the
solution is uniquely determined by its largest coordinate) fails in
our case, more precisely:
\begin{thm}
  \label{t:fixmax} For $c \in \cC_3$ and $c > 1$, the number of fair
  games with $c$ as the largest coordinate is $2^{m-1}$ where $m$ is
  the number of distinct prime factors of $2c+1$ other than 3.
\end{thm}
\begin{proof}
  The coordinates of nontrivial 3-color fair games are distinct, so it
  follows from~\eqref{eq:3para} that $c$ is the largest coordinates if
  and only if
  \[
  c \ge \frac{u^2 + (2c+1)|u| + c(c-1)}{2(2c+1)} +1.
  \]
  A simple calculation shows that the above inequality is equivalent
  to $|u| \le c$. Hence, proving the theorem boils down to counting
  the number of solutions to the congruence $u^2 \equiv -c(c-1) \mod
  (2c+1)$, or equivalently $u^2 \equiv -3 \mod (2c+1)$ in a complete
  set of representatives ($-c \le u \le c$). By
  Theorem~\ref{t:partition}, the prime factorization of $2c+1$ is of
  the form $3^{v_0}p_1^{v_1}\cdots p_m^{v_m}$ where $v_0 =0$ or $1$
  and $p_j \equiv 1 \mod 3$ ($1 \le j \le m$). Note that $m \ge 1$
  since $c > 1$. The Chinese reminder theorem yields a 1-to-1
  correspondence between solutions of $u^2 \equiv -3 \mod (2c+1)$ and
  the solutions of the system $u^2 \equiv -3 \mod{p_j^{v_j}}$ ($1 \le
  j \le m$) together with the congruence $u^2 \equiv -3 \equiv 0 \mod
  3$ depending on whether $v=1$ or not. But last congruence has only
  one solution hence its presence will not affect the total number of
  solutions which is the product of the number of solutions of each
  congruence in the system. Since $p_j \neq 3$ ($1 \le j \le m$), the
  number of solutions for $u^2 \equiv -3 \mod p_j^{v_j}$ is the same
  as that for $u^2 \equiv -3 \mod p_j$ which is precisely two since
  $p_j \equiv 1 \mod 3$. Therefore, we conclude that the system has
  exactly $2^m$ solutions. Finally, it is clear from~\eqref{eq:3para}
  that $\pm u$ give rise to the same fair game except with the first
  two coordinates permuted. Thus, up to permutation of coordinates,
  there are $2^{m}/2 = 2^{m-1}$ fair games with $c$ as the largest
  coordinate.
\end{proof}
We conclude this section by computing the natural density of $\cC_3$.
The {\em natural density} of a set of natural numbers $A$ is defined
to be $\lim_{k \to \infty} |A(k)|/k$ whenever the limit
exists. Statement~\eqref{i:C2} of Corollary~\ref{c:asymn=2} states
that $\cC_2$ has density zero. The description of $\cC_3$ in
Theorem~\ref{t:partition} allows us to show that the same phenomenon
occurs in the case $n=3$ as well.

\begin{thm}
  \label{t:density}
  The natural density of $\cC_3$ is zero.
\end{thm}
\begin{proof}
  The map $a \mapsto 2a+1$ is a bijection between $\cC_3$ and
  $P_{1}^{\ge 0} \cup 3P_{1}^{\ge 0}$. Therefore, the density of
  $\cC_3$ is twice the density of $P_{1}^{\ge 0} \cup 3P_{1}^{\ge 0}$
  if the latter exists. Since $P_{1}^{\ge 0}$ and $3P_{1}^{\ge 0}$ are
  disjoint and the density of $3P_{1}^{\ge0}$, if it exists, is
  one-third that of $P_1^{\ge 0}$, it suffices to show that
  $P_{1}^{\ge 0}$ has density $0$.  Applying Proposition~9.64 and
  Lemma~11.8 in~\cite{ntd} to the set $P_1^{\ge0} = P_1^{\ge 0}P_0^0
  P_{-1}^0$, we see that it is enough to show that the series $\sum_{p
    \in P_{-1}} 1/p$ diverges. But this assertion is an immediate
  consequence of Dirichlet Theorem of primes in arithmetic
  progressions~\cite[Chapter~VI Theorem 2]{ser}) which, in particular,
  asserts that
  \[
  \lim_{s \to 1^+} \frac{\sum_{p \in P_{-1}} 1/p^s}{\log (1/(s-1))} =
  \frac{1}{\phi(3)} = \frac{1}{2}.
  \]
\end{proof}

\section{Asymptotic Behavior}
\label{s:asy}
Comparing to the binary case, determining the asymptotic behavior of
$|\cF_n(k)|$ ($n \ge 3$) seems to be much harder problem. In fact, we
will only show that $|\cF_3(k)|$ is $\Theta(k)$. Our strategy is to
relate the equation for 3-color fair games to the Lorentzian form
$L_3(w) = w_1^2 + w_2^2 - w_3^2$ then apply the results
from~\cite{lorentzian} (for more information on distribution of
integral points on affine homogenous varieties, see~\cite{drs}
and~\cite{mb}). Let $W$ be the set of integral solutions of $L_3(w) =
-3$.
\begin{lem}
  \label{l:lorentz}
  For any $(w_1, w_2, w_3) \in \sW$, we have
  \begin{enumerate}
  \item \label{i:abs} $|w_3| > |w_1|, |w_2|$.
    
  \item \label{i:oddcoord} exactly one of the $w_i$ is odd; moreover,
    it must be either $w_1$ or $w_2$.
    
  \item \label{i:w_2odd} If $w_2$ is the odd coordinate that $w_1+w_3
    \equiv w_1-w_3 \equiv 2 \mod 4$.

  \item \label{i:neq} $w_1 \neq w_2$.
  \end{enumerate}
\end{lem}
\begin{proof}
  The first statement is immediate. The second statement is clear by
  arguing mod 8. Since $w_2$ is odd, $w_1^2 -w_3^2 \equiv 0 \mod
  4$. Since $w_3$ is even, $w_1 + w_3 \equiv w_1 - w_3 \mod 4$. They
  must be both congruent to $2$ mod $4$; otherwise $w_1^2 - w_3^2
  \equiv 0 \mod 8$ making $-3$ a square mod $8$, contradiction. This
  establishes the third statement. The last statement is true since
  $2w_1^2 -w_3^2 =-3$ is not solvable mod 3.
\end{proof}
Let $\approx$ be the equivalence relation on $\sW$ identifying the
elements $(w_1, w_2, w_3)$ and $(w_2, w_1, w_3)$. It follows from
Lemma~\ref{l:lorentz}~\eqref{i:neq} that the canonical map from $\sW$
to $\sW/\!\approx$ is 2-to-1. If we identify the equivalence classes
with those elements of $\sW$ with an odd second coordinate, then
\begin{prop}
  \label{p:equiv}
  The map given by
    \begin{equation}
      \label{eq:x-to-w}
      w_1 = 2(x_2 - x_3), \quad w_2 = 2(x_1 - x_2 -x_3) -1, \quad w_3
      = 2(x_2 + x_3 +1)
    \end{equation}
    is a 1-to-1 correspondence between $\cS_3$ and $\sW/\!\approx$.
    Moreover, if the coordinates of the elements of $\cS_3$ are listed
    in ascending order, then elements of $\cF_3$ correspond to those
    elements of $\sW/\!\approx$ with $w_1, w_2 \le 0$ and $w_3 \ge 0$.
\end{prop}
\begin{proof}
  The rational inverse of the map in~\eqref{eq:x-to-w} is given by  
  \begin{equation}
    \label{eq:w-to-x}
    x_1 = \frac{w_2 + w_3 -1}{2}, \quad x_2 = \frac{w_1 + w_3 -2}{4},
    \quad x_3 = \frac{w_3 - w_1 -2}{4}.
  \end{equation}
  By Lemma~\ref{l:lorentz}~\eqref{i:oddcoord} and~\eqref{i:w_2odd}, it
  actually preserve integral points. This establishes the map
  in~\eqref{eq:x-to-w} is a 1-to-1 correspondence between $\cS_3$ and
  $\sW/\!\approx$. Moreover, the images of elements of $\cF_3$ (as
  ascending triples) under~\eqref{eq:x-to-w} clearly satisfy $w_1, w_2
  \le 0$ and $w_3 \ge 0$. Conversely, by
  Lemma~\ref{l:lorentz}~\eqref{i:abs}, $|w_3| \ge |w_2| + 1$ and
  $|w_3| \ge |w_1| + 2$ since both $w_1$ and $w_3$ are odd
  (Lemma~\ref{l:lorentz}~\eqref{i:oddcoord}). Therefore, if
  $(w_1,w_2,w_3) \in \sW$ with $w_1, w_2 \le 0$ and $w_3 \ge 0$ then
  the corresponding $(x_1,x_2,x_3)$ is in $\cF_3$.
\end{proof}

\begin{thm}
  \label{t:Theta}
  There exists positive constants $c_1,c_2$ such that $ c_1k \le
  |\cS_3(k)| \le c_2k$ for $k$ sufficiently large, i.e. $|\cS_3(k)| =
  \Theta(k)$. Similarly, $|\cF_3(k)| =\Theta(k)$.
\end{thm}
\begin{proof}
  The idea is simple: the sphere of radius $k$ centered at the origin
  maps to an ellipsoid (centered at $(0,-1,2)$)
  under~\eqref{eq:x-to-w}. For $k$ sufficiently large, it is enveloped
  between spheres centered at the origin. Clearly, the radii of these
  spheres can be chosen as linear functions of the length of the axes
  of the ellipsoid which are in turn linear in $k$ since the
  transformation given in~\eqref{eq:x-to-w} is affine. So by
  Proposition~\ref{p:equiv}, the proof is complete once we show that
  $|\sW(k)|$ is asymptotic to a linear function in $k$. And this last
  statement follows from Formula~(3) in~\cite{lorentzian} which
  asserts that $|\sW_3(k)| \sim (4\sqrt{6}/3)k$.
\end{proof}

\begin{rmk}
\label{r:asy} \hfill{}
\begin{enumeratei}
\item Since the maps and the equations are all explicit, one can
  provide the constants $c_1$ and $c_2$ in Theorem~\ref{t:Theta}
  explicitly. However, we will leave the computation for the
  interested readers.

\item In a sense, one gets a cleaner result without extra efforts if
  one is satisfied by counting the number of solutions inside the
  ellipsoids that are the images of spheres under the map given
  in~\eqref{eq:w-to-x}. To be more precise, let $\cS_3'(k)$ and
  $\cF_3'(k)$ be the set elements of $\cS_3$ and $\cF_3$ inside the
  image of the sphere of radius $k$ centered at the origin under the
  map given in~\eqref{eq:w-to-x}. Then again using Formula~(3),
  Table~(1) in~\cite{lorentzian} and Proposition~\ref{p:equiv}, one
  gets
  \[
  |\cS_3'(k)| = \frac{1}{2}|\sW(k)| \sim \frac{2\sqrt{6}}{3}k, \quad
  |\cF'_3(k)| \sim \frac{6}{8}\frac{1}{2}|\sW(k)| =
  \frac{\sqrt{6}}{2}k.
  \]

\item Here is how we arrive to the map in~\eqref{eq:x-to-w}: There is
  a general method of solving quadratic Diophantine equations given by
  Grunewald and Segal in~\cite{gs} and~\cite{sa}\footnote{However,
    their method is not uniform in the number of variables.}. The
  first step transforms the fair game equation into the equivalent
  system
  \begin{equation}
  \label{eq:n-system}
  Q_n(\bds{z})=-n(n-2), \quad z_i \equiv 1 \mod 2(n-2).
  \qquad (1 \le i \le n)
  \end{equation}
  where $Q_n$ is the quadratic form in $n$-variables with diagonal
  entries 1 and off-diagonal entries -1. When $n=3$, one checks
  readily that congruences in~\eqref{eq:n-system} are implies by
  $Q_3(\bds{z}) = -3$. And the Lorentzian form $L_3(\bds{w})$ is
  obtained by diagonalizing $Q_3$.  Taking the composition of these
  transformations yields the map in~\eqref{eq:x-to-w}.

\item In the ternary case, the description of the solution set given
  by Grunewald and Segal's method relates quite beautifully to
  ours. We encourage the reader to pursuit their original papers
  (see~\cite{cas} for the necessary backgrounds). Just to give a
  little enticement, let us remark that $\cS_3$ will correspond to a
  single orbit under the integral orthogonal group of a suitable
  quadratic form. While $\cS_3^{+}$ (i.e. $\cF_3$) and $\cS_3^{-}$
  will correspond to two orbits of a subgroup of the orthogonal group.

\item Using the same idea to study $|\cF_n(k)|$ for $n \ge 4$ becomes
  more problematic. First, it is not clear to us how to take into
  account of the congruences in~\eqref{eq:n-system}. Moreover, even
  though $Q_n$ and the Lorentzian form in $n$-variables have the same
  signature, namely $n-2$, the transformation taking one to the other
  in general does not preserve integral points.
\end{enumeratei}
\end{rmk}

\section{Odds and Ends} 
\label{s:misc} 
The last section is dedicated to various results that do not quite fit
in previous sections.
\begin{prop}
  \label{p:mod4}
  The sum of coordinates of any element of $S_n$ is congruent to
  either 0 or 1 mod 4.
\end{prop}
\begin{proof}
  Equation~\eqref{eq:fair} is simply $s^2 -s = 4p$.
\end{proof}
\begin{prop}
  \label{p:mod3}
  Every vertex in the connected component of $\bds{0}$ in $\cS_{n}$ is
  congruent mod 3 to either $\bds{0}$ or
  $\bds{e}_j:=(\tuple{0}{1,\dots,0})$ for some $1 \le j \le n$.
\end{prop}
\begin{proof}
  The proposition is clearly true for $\bds{0}$. Now suppose it is true
  for every vertex of distance $m$ from $\bds{0}$. Let ${\bds a}'=(a_1,
  \ldots, a_i', \ldots, a_n)$ be a vertex of distance $m+1$ from
  $\bds{0}$ and is adjacent to $\bds{a}=(a_1, \ldots, a_i, \ldots, a_n)$
  which is of distance $m$ from $\bds{0}$. By the induction hypothesis,
  $\bds{a}$ is congruent to either $\bds{0}$ or $\bds{e}_j$ mod $3$ for
  some $1 \le j \le n$. Since $a_i + a_i'=1+ 2s_i(\bds{a}) \equiv
  1-s_i(\bds{a}) \mod 3$, we have the following 3 cases
  \begin{enumerate}
  \item ${\bds a}' \equiv \bds{e}_i \mod 3$ if $\bds{a} \equiv \bds{0}
    \mod 3$,
  \item ${\bds a}' \equiv \bds{e}_j \mod 3$ if $\bds{a} \equiv \bds{e}_j
    \mod 3$ and $i \neq j$, or
  \item ${\bds a}' \equiv \bds{0} \mod 3$ if $\bds{a} \equiv \bds{e}_j
    \mod 3$ and $i = j$.
  \end{enumerate}
  This establishes the proposition by induction.
\end{proof}
Recall that for $n \ge 4$ any natural number, in particular 2, can be a
coordinate of an $n$-color fair game (Proposition~\ref{p:C_n}). Thus,
Proposition~\ref{p:mod3} gives another proof of the fact that $\cF_n$ is
disconnected for $n \ge 4$.

The nontrivial 3-color fair games form a full binary tree with the
nontrivial $2$-color fair games embeds as a branch. It is easy to see
that among the nodes of distance $k$ from the root $(0,1,3)$, the one
with the smallest norm is $\big(0, \binom{k+2}{2},
\binom{k+3}{2}\big)$ and the one with the largest norm is $(m_k,
m_{k+1}, m_{k+2})$ where $(m_i)$ is the sequence defined recursively
by $(m_0,m_1,m_2) = (0,1,3)$ and $m_{i+3}= 2(m_{i+1} + m_{i+2})+1 -
m_i$ for $i \ge 0$. It turns out that the $m_i$'s have an intimate
relation with the Fibonacci numbers.
\begin{prop}
  \label{p:fibonacci}
  Let $f_i$ be the $i$-th Fibonacci number, then for any $k \ge 0$
  \begin{equation*}
    m_k =
    \begin{cases}
      f_k^2 &\text{if $k$ is odd} \\
      f_k^2 -1 &\text{if $k$ is even}
    \end{cases}
  \end{equation*}
\end{prop}
\begin{proof}
  The Fibonacci numbers are defined inductively by $f_0=f_1=1$, and
  $f_{i+2}=f_{i} + f_{i+1}$ ($i \ge 0$). Note that
  \begin{align*}
    2f_{i+1}f_{i+2} &= f_{i+1}(f_{i+2} + f_{i} + f_{i+1}) =
    f_{i+1}(f_{i+2} +f_{i}) + f_{i+1}^2 \\
    &=(f_{i+2}-f_{i})(f_{i+2}+f_{i}) +f_{i+1}^2 \\
    &= f_{i+2}^2 + f_{i+1}^2 - f_{i}^2.
  \end{align*}
  Therefore,
  \begin{align*}
    f_{i+3}^2 &=(f_{i+1} + f_{i+2})^2 = f_{i+1}^2 + 2f_{i+1}f_{i+2} +
    f_{i+2}^2 \\
    &= 2f_{i+1}^2 + 2f_{i+2}^2 - f_{i}^2.
  \end{align*}
  The proposition now follows from an easy induction. The base cases
  are immediate. Suppose $k \ge 3$ and the proposition is true for all
  $0 \le i < k$.  When $k$ is odd, we have
  \begin{align*}
    m_{k} &= 2(m_{k-2} + m_{k-1}) +1 - m_{k-3} \\
    &= 2(f_{k-2}^2 + f_{k-1}^2 -1) +1 -(f_{k-3}^2 -1) \\
    &= 2f_{k-2}^2 + 2f_{k-1}^2 -f_{k-3}^2 \\
    &= f_{k}^2.
  \end{align*}
  A similar computation shows that $m_k = f_k^2 -1$ when $k$ is even.
\end{proof}

We end the article with another curious ``by-product'' of our results.
\begin{prop}
  \label{p:factors}
  For $m \ge 0$,
  \begin{enumerate}[\upshape{(}1\upshape{)}]
  \item \label{i:c2-factor} $m^2 + m +1 \in P_1^{\ge 0} \cup 3P_1^{\ge
      0}$.
    
  \item \label{i:fi-factor} $2f_m^2 - (-1)^m \in P_1^{\ge 0} \cup
    3P_1^{\ge 0}$.
  \end{enumerate}
\end{prop}
\begin{proof}
  Triangular numbers $m(m+1)/2$ ($m \ge 0$) appear as coordinates of the
  2-color (Theorem~\ref{t:n=2}) and hence 3-color fair games. By
  Theorem~\ref{t:partition}, we have $m^2 +m +1 = 2(m(m+1)/2) +1 \in
  P_1^{\ge 0} \cup 3P_1^{\ge 0}$. Similarly,
  Statement~\eqref{i:fi-factor} follows from
  Proposition~\ref{p:fibonacci} and Theorem~\ref{t:partition}.
\end{proof}

\end{document}